\documentclass[12pt]{amsart}
\newtheorem{theorem}{Theorem}
\newtheorem{proposition}[theorem]{Proposition}
\newtheorem{lemma}[theorem]{Lemma}

\def\ds{\displaystyle}

\title[The Carath\'eodory metric on the symmetrized polydisc]
{Estimates of the Carath\'eodory metric on the symmetrized
polydisc}

\author[N.~Nikolov, P.~Pflug, P.~J.~Thomas and W.~Zwonek]
{Nikolai Nikolov, Peter Pflug, Pascal J. Thomas and W\l odzimierz
Zwonek}

\address
{Institute of Mathematics and Informatics\\ Bulgarian Academy of
Sciences\\ Acad. G. Bonchev 8, 1113 Sofia,
Bulgaria}\email{nik@math.bas.bg}

\address{Carl von Ossietzky Universit\"at Oldenburg\\
Fachbereich Mathe\-ma\-tik\\ Postfach 2503\\ D-26111 Oldenburg,
Germany}\email{pflug@mathematik.uni-oldenburg.de}

\address{Laboratoire Emile Picard, UMR CNRS 5580\\
Universit\'e Paul Sabatier, 118 Route de Narbonne\\ F-31062
Toulouse Cedex, France} \email{pthomas@cict.fr}

\address{Instytut Matematyki, Uniwersytet Jagiello\'nski, Reymonta 4,
30-059 Krak\'ow, Poland}\email{Wlodzimierz.Zwonek@im.uj.edu.pl}



\begin{document}

\begin{thanks} {The first and second named authors were supported by grants from
DFG (DFG-Projekt 227/8); the fourth-named author was supported by
the KBN Research Grant No. 1 PO3A 005 28 and Alexander von
Humboldt Foundation.}
\end{thanks}

\begin{abstract} Estimates for the Carath\'eodory metric on the symmetrized
polydisc are obtained. It is also shown that the Carath\'eo\-dory
and Kobayashi distances of the symmetrized three-disc do not
coincide.
\end{abstract}

\maketitle

\section{Introduction}

A consequence of the fundamental Lempert theorem (see \cite{Lem})
is the fact that the Carath\'eodory distance and the Lempert
function coincide on any domain $D\subset\Bbb C^n$ with the
following property $(\ast)$ (cf. \cite{Jar-Pfl}):

\noindent ($\ast$) {\it  $D$ can be exhausted by domains which are
biholomorphic to convex domains.}

For more than 20 years it has been an open question whether the
converse of the above result is true in some reasonable class of
domains (e.g. in the class of bounded pseudoconvex domains). In
other words, does the equality between the Carath\'eodory distance
$c_D$ and the Lempert function $\tilde k_D$ of a bounded
pseudoconvex domain $D$ imply that $D$ satisfies property $(\ast)$.

The first counterexample, the so-called symmetrized bidisc $\Bbb
G_2,$ has been recently discovered and discussed in a series of papers
(see \cite{Agl-You}, \cite{Cos1}, \cite{Cos2} and \cite{Edi}, see
also \cite{Jar-Pfl}).

In fact, it was proved that $c_{\Bbb G_2}$ and $\tilde k_{\Bbb G_2}$
coincide with a natural distance $p_{\Bbb G_2}$ related to (the geometry of)
$\Bbb G_2$.

The symmetrized polydisc $\Bbb G_n$ ($n\ge 3$) can also be endowed
with a similar distance $p_{\Bbb G_n}$ which does not exceed
$c_{\Bbb G_n}.$ Using $p_{\Bbb G_n},$ three of the authors have
recently shown that $\tilde k_{\Bbb G_n}$ is not a distance (see
\cite{NPZ}); in particular, $\Bbb G_n$ does not satisfy property
$(\ast)$ (for a direct proof of this fact see \cite{Nik}). They have
also proved that the Kobayashi distance of $\Bbb G_n$ does not
coincide with $p_{\Bbb G_n}.$

In the present paper we improve this result showing that $c_{\Bbb
G_n}(0;\cdot)\neq p_{\Bbb G_n}(0;\cdot).$ The proof is based on
the comparison of the infinitesimal version of these distances at
the origin, $\gamma_{\Bbb G_n}(0;\cdot)$ and $\rho_n$, where
$\gamma_{\Bbb G_n}$ is the Carath\'eodory-Reiffen metric of $\Bbb
G_n.$ We also give lower and upper bounds for $\gamma_{\Bbb
G_{2n+1}}(0;e_2)$ (where $e_2$ is the second basis vector). The
bounds give an asymptotic estimate for $\gamma_{\Bbb
G_{2n+1}}(0;e_2)$ with an error of the form $o(n^{-3}).$ Finally,
estimating more precisely the value of $\gamma_{\Bbb G_{3}}$ at
the point $(0;e_2)\in\Bbb G_3\times\Bbb C^3,$ we obtain that it is
smaller than the infinitesimal version of the Kobayashi distance
at the same point which implies that the Kobayashi distance does
not coincide with the Carath\'eodory distance on $\Bbb G_3.$
\smallskip

\noindent{\bf Acknowledgements.} We thank Dr.~Pencho Marinov for
the computer programmes helping us to obtain the estimates in the
last section of the paper.

\section{Background}

Let $\Bbb D$ be the unit disc in $\Bbb C.$ Let
$\sigma_n=(\sigma_{n,1},\dots,\sigma_{n,n}):\Bbb C^n\to\Bbb C^n$
be defined as follows: $$\sigma_{n,k}(z_1,\dots,z_n)=\sum_{1\le
j_1<\dots<j_k\le n}z_{j_1}\dots z_{j_k},\quad 1\le k\le n.$$ The
domain $\Bbb G_n=\sigma_n(\Bbb D^n)$ is called {\it the
symmetrized $n$-disc}.

Recall now the definitions of the Carath\'eodory pseudodistance,
the Carath\'eodory-Reiffen  pseudometric, the Lempert function and
the Ko\-bayashi-Royden pseudometric of a domain $D\subset\Bbb C^n$
(cf. \cite{Jar-Pfl}): $$ \aligned
c_D(z,w)&:=\sup\{\tanh^{-1}|f(w)|:f\in\mathcal O(D,\Bbb
D),f(z)=0\},\\ \gamma_D(z;X)&:=\sup\{|f'(z)X|:f\in\mathcal
O(D,\Bbb D),f(z)=0\},\\ \tilde
k_D(z,w)&:=\inf\{\tanh^{-1}|\alpha|:\exists\varphi\in\mathcal
O(\Bbb D,D): \varphi(0)=z,\varphi(\alpha)=w\},\\
\kappa_D(z;X)&:=\inf\{\alpha\ge 0:\exists\varphi\in\mathcal O(\Bbb
D,D): \varphi(0)=z,\alpha\varphi'(0)=X\},
\endaligned
$$ where $z,w\in D,$ $X\in\Bbb C^n.$ The Kobayashi pseudodistance
$k_D$ (respectively, the Kobayashi--Buseman pseudometric
$\hat\kappa_D$) is the largest pseudodistance (respectively,
pseudonorm) which does not exceed $\tilde k_D$ (respectively,
$\kappa_D$).

It is well-know that $c_D\le k_D\le\tilde k_D,$
$\gamma_D\le\hat\kappa_D\le\kappa_D,$ and
$$\gamma_D(z;X)=\lim_{\Bbb C_\ast\ni t\to 0}\frac{c_D(z,z+tX)}{t}\
\hbox{(cf. \cite{Jar-Pfl})}.$$ Moreover, if $D$ is taut, then
$$\ds\kappa_D(z;X)= \lim_{\Bbb C_\ast\ni t\to 0}\frac{\tilde
k_D(z,z+tX)}{t}\ \hbox{(see \cite{Pang}},$$
$$\hat\kappa_D(z;X)=\lim_{\Bbb C_\ast\ni t\to
0}\frac{k_D(z,z+tX)}{t}\ \hbox{(see \cite{KobM})}.$$

Note that $\Bbb G_n$ is a hyperconvex domain (see \cite{Edi-Zwo})
and, therefore, a taut domain.

In the proofs below we shall need some mappings defined on $\Bbb
G_n.$

For $\lambda\in\overline{\Bbb D}$, $n\ge 2$, one may define the
rational mapping $p_{n,\lambda}$ as follows $$ p_{n,\lambda}(z):=(\tilde
z_1(\lambda),\ldots,\tilde z_{n-1}(\lambda))=\tilde z(\lambda)
\in\Bbb C^{n-1},\quad z\in\Bbb C^n, n+\lambda z_1\neq 0,$$ where
$\ds\tilde
z_j(\lambda)=\frac{(n-j)z_j+\lambda(j+1)z_{j+1}}{n+\lambda z_1},$
$1\le j\le n-1$. Then $z\in\Bbb G_n$ if and only if $\tilde
z(\lambda)\in\Bbb G_{n-1}$, $n+\lambda z_1\neq 0$ for any $\lambda\in\overline{\Bbb D}$
(see Corollary 3.4 in \cite{Cos3}).

We may also define for
$\lambda_1,\ldots,\lambda_{n-1}\in\overline{\Bbb D}$ the rational
function $$ f_{\lambda_1,\ldots,\lambda_{n-1}}=
p_{2,\lambda_1}\circ\ldots\circ p_{n,\lambda_{n-1}}.$$ Observe
that
$$f_\lambda(z):=f_{\lambda,\dots,\lambda}(z)=\frac{\sum_{j=1}^{n}jz_j
\lambda^{j-1}}{n+\sum_{j=1}^{n-1}(n-j)z_j\lambda^j}.$$ By Theorem
3.2 in \cite{Cos3}, $z\in\Bbb G_n$ if and only if
$\ds\sup_{\lambda\in\overline{\Bbb D}}|f_\lambda(z)|<1.$ In fact,
by Theorem 3.5 in \cite{Cos3}, if $z\in\Bbb G_n,$ then the last
supremum is equal to
$\ds\sup_{\lambda_1,\dots,\lambda_{n-1}\in\overline{\Bbb D}
}|f_{\lambda_1,\dots,\lambda_{n-1}}(z)|$.

It follows that $$c_{\Bbb G_n}(z,w)\ge p_{\Bbb G_n}(z,w):=
\max_{\lambda_1,\dots,\lambda_{n-1}\in\Bbb T}|p_{\Bbb
D}(f_{\lambda_1,\dots,\lambda_{n-1}}(z),f_{\lambda_1,\dots,\lambda_{n-1}}(w))|,$$
where $\Bbb T=\partial\Bbb D$ and $p_{\Bbb D}$ is the Poincar\'e
distance. Observe that $p_{\Bbb G_n}$ is a distance on $\Bbb G_n$.

Let $e_1,\dots,e_n$ be the standard basis of $\Bbb C^n$ and $\ds
X=\sum_{j=1}^n X_je_j.$ Set $$\tilde
f_{\lambda}(X)=\frac{\sum_{j=1}^{n}jX_j \lambda^{j-1}}{n}\hbox{\ \
and\ \ }\rho_n(X):=\max_{\lambda\in\Bbb T}|\tilde f_\lambda(X)|.$$
Then the last inequality above implies that $$\gamma_{\Bbb
G_n}(0;X)\ge\lim_{\Bbb C_\ast\ni t\to 0}\frac{p_{\Bbb
G_n}(0,tX)}{|t|}=\rho_n(X).$$

Let $L_{k,l}$  be the span of $e_k$ and $e_l.$ Note that if $X \in
L_{k,l}$, $k\neq l$, then $$\rho_n(X)=\frac{k|X_k|+l|X_l|}{n}.$$

For $n=2$ one has equalities $k_{\Bbb G_2}=c_{\Bbb G_2}=p_{\Bbb G_2}$
(see \cite{Agl-You}, \cite{Cos1}). On the other hand, we have the
following (see \cite{NPZ}).

\begin{proposition} Let $n\ge 3.$

(a) If $k$ divides $n,$ then $\ds\kappa_{\Bbb
G_n}(0;e_k)=\rho_n(e_k).$ Therefore, if $l$ also divides $n,$ then
$\hat\kappa_{\Bbb G_n}(0;X)=c_{\Bbb G_n}(0;X)=\rho_n(X)$ for any
$X\in L_{k,l}.$

(b) If $X\in L_{1,n}\setminus(L_{1,1}\cup L_{n,n}),$ then
$\kappa_{\Bbb G_n}(0;X)>\rho_n(X).$

(c) If $k$ does not divide $n,$ then $\hat\kappa_{\Bbb
G_n}(0;e_k)>\rho_n(e_k).$

In particular, $\tilde k_{\Bbb G_n}(0,\cdot)\neq k_{\Bbb
G_n}(0,\cdot)$ and $k_{\Bbb G_n}(0,\cdot)\neq p_{\Bbb
G_n}(0,\cdot).$
\end{proposition}

In the next section we shall prove a stronger inequality than that
in Proposition 1 (c).

\section{If $k$ does not divide $n,$ then $\gamma_{\Bbb G_n}(0;e_k)>\rho_{\Bbb G_n}(e_k)$}
Our first aim is the proof of a result, which implies the inequality between $c_{\Bbb G_n}$ and $p_{\Bbb G_n}$, $n\geq 3$.

\begin{proposition} If $k$ does not divide $n\ge 3,$ then
$\gamma_{\Bbb G_n}(0;e_k)>\rho_n(e_k).$ In particular, $c_{\Bbb
G_n}(0,\cdot)\neq p_{\Bbb G_n}(0,\cdot).$
\end{proposition}

\begin{proof}
Let $\root k\of 1=\{\xi_1,\dots,\xi_k\}.$ For $z\in\overline{\Bbb
G_n}$ and $\lambda\in\overline{\Bbb D},$ such that the denominator
in the formula below does not vanish, set $$g_z(\lambda):=\lambda
f_\lambda(z)=\frac{\sum_{j=1}^{n}jz_j
\lambda^j}{n+\sum_{j=1}^{n-1}(n-j)z_j\lambda^j}$$ and
$$g_{z,k}(\lambda)=\frac{\sum_{j=1}^kg_z(\xi_j\lambda)}{k\lambda^k}.$$
The equalities $\sum_{j=1}^k\xi_j^m=0,$ $m=1,\dots,k-1,$ and the
Taylor expansion of $g_{z,k}$ show that this function can be
extended at $0$ as $g_{z,k}(0)=P_k(z),$ where $P_k$ is a
polynomial with $\ds\frac{\partial P_k}{\partial
z_k}_{|_{z=0}}=\frac{k}{n}$ and
$$P_k(tw_1,t^2w_2,\dots,t^nw_n)=t^kP(w),\
w=(w_1,w_2,\dots,w_n)\in\Bbb C^n,t\in\Bbb C.$$ It follows by the
maximum principle that $g_{z,k}\in\mathcal O(\Bbb D,\overline{\Bbb
D})$. In particular, $|P_{k}(z)|\le 1.$ To prove the desired
inequality, it is enough to show that $|P_k(z)|<1$ for any
$z\in\overline{\Bbb G_n}.$ Assume the contrary, that is,
$P_k(z)=e^{i\theta}$ for some $\theta\in\Bbb R$ and some
$z\in\overline{\Bbb G_n}.$ Then the maximum principle and the
triangle inequality implies that $g_z(\xi_j\lambda)=
e^{i\theta}\lambda^k,$ $\lambda\in\Bbb T,$ $1\le j\le k.$ In
particular, $g_z(\lambda)=e^{i\theta}\lambda^k,$ that is
$$\sum_{j=1}^{n}jz_j\lambda^j=
e^{i\theta}(n\lambda^{k}+\sum_{j=1}^{n-1}(n-j)z_j\lambda^{k+j}).$$
Comparing the corresponding coefficients of these two polynomials
of $\lambda,$ we get that $\ds z_k=e^{i\theta}\frac{n}{k},$
$z_{n+1-k}=\dots=z_{n-1}=0$ and
$$(k+j)z_{k+j}=e^{i\theta}(n-j)z_j,\ 1\le j\le n-k.$$ The last
relations imply that $\ds z_{kl}=e^{i\theta}\binom{n/k}{l},$ $1\le
l\le[n/k].$ On the other hand, since $k$ does not divide $n,$ then
$n-k<k[n/k]<n$ and hence $z_{k[n/k]}=0$ -- a contradiction.
\end{proof}

\noindent{\bf Remarks.} It will be interesting to know whether
$\hat\kappa_{\Bbb G_n}(0;\cdot)\neq\gamma_{\Bbb G_n}(0;\cdot)$ and
hence $k_{\Bbb G_n}(0,\cdot)\neq c_{\Bbb G_n}(0,\cdot)$ for any
$n\ge 4.$ In the last section we shall prove these inequalities
for $n=3.$

\section{Estimates for $\ds\gamma_{\Bbb G_{2n+1}}(0;e_2)$}

Let $n$ and $k$ be positive integers, $k\le n.$ Note that
$$\kappa_{\Bbb G_n}(0;e_k)\le\kappa_{\Bbb
G_{k[n/k]}}(0;e_k)=\frac{1}{[n/k]}.$$ Thus,
$$\frac{k}{n}\le\gamma_{\Bbb G_n}(0;e_k)\le\kappa_{\Bbb
G_n}(0;e_k)\le\frac{1}{[n/k]}.$$ Therefore, one has that
$$\lim_{n\to\infty}n\gamma_{\Bbb
G_n}(0;e_k)=\lim_{n\to\infty}n\kappa_{\Bbb G_n}(0;e_k)=k.$$

Let now $n\ge 3$ be odd. It follows that $$\gamma_{\Bbb
G_n}(0;e_2)\le\kappa_{\Bbb G_n}(0;e_2)\le\frac{2}{n-1}.$$ On the
other hand, $\ds\frac{2}{n}<\ds\gamma_{\Bbb G_n}(0;e_2)$ by
Proposition 1. The aim of this section is to improve both
estimates.

To obtain a more precise upper bound, we shall need the following
definition. Let $k_1\le\dots\le k_n$ be positive integers. For
$\lambda\in\Bbb C,$ define the mapping $$\pi_\lambda:\Bbb
C^n\ni(z_1,\dots,z_n)\to
(\lambda^{k_1}z_1,\dots,\lambda^{k_n}z_n)\in\Bbb C^n.$$ We shall
say that a domain $D\subset\Bbb C^n$ is $(k_1,\dots,k_n)$-balanced
if $\pi_\lambda(z)\in D$ for any $\lambda\in\overline{\Bbb D}$ and
any $z\in D.$ For such a domain $D$ and any $j=1,\dots,n,$ denote
by $\mathcal P_j$ the set of polynomials $P$ with $\sup_D|P|\le
1$ and $P\circ\pi_\lambda=\lambda^{k_j}P,$ and by $\mathcal L_j$
the span of the vectors $e_j,\dots,e_l,$ where $l\ge j$ is the
maximal integer with $k_l=k_j.$ The proof of Proposition 1 implies
the following result.

\begin{proposition} If $D\subset\Bbb C^n$ is a
$(k_1,\dots,k_n)$-balanced domain and $X\in\mathcal L_j,$ $1\le
j\le n$, then $\gamma_D(0;X)=\sup\{|P'(0)X|:P\in\mathcal P_j\}.$
\end{proposition}

\noindent{\bf Remarks.} (i) One can obtain a similar description
for any Reiffen pseudometric of higher order (for the definition
see the next section).

(ii) A consequence of Proposition 3 is the well-know fact that if
$D$ is a balanced domain, that is, $k_1=\dots=k_n=1,$ then
$\gamma_D(0;X)=\hat h_D(X),$ where $\hat h_D$ is the Minkowski
function of the convex hull of $D.$

(iii) Another consequence of Proposition 3 is the formula
$$\gamma_{\Bbb G_n}(0;e_2)=\frac{1}{\inf_{c\in\Bbb
C}\max_{z\in\partial G_n}|z_2+cz_1^2|}.\leqno{(1)}$$ Despite of
(1), it is difficult to find explicitly $\gamma_{\Bbb G_n}(0;e_2)$
for odd $n\ge 3$ (see the last section).

(iv) Note that in the case of an even $n$ the extremal polynomials
for $\ds\gamma_{\Bbb G_n}(0;e_2)=\frac{2}{n}$ are not unique up to
a rotation. Namely, the proof of Proposition 2 delivers the
polynomial $\ds\frac{2}{n}z_2-\frac{n-1}{n^2}z_1^2,$ but
$\ds\frac{2}{n}z_2-\frac{1}{n}z_1^2$ is also an extremal
polynomial.

\begin{proposition} If $n\ge 3$ is odd, then
$$\frac{2}{n}\left(1+\frac{2}{(n-1)(n+2)}\right)<\gamma_{\Bbb
G_n}(0;e_2)<\frac{2}{n}\left(1+\frac{2}{(n-1)(n+1)}\right).$$
\end{proposition}

\begin{proof} {\it The lower bound:}
First, we shall see that for the polynomial $\ds P_n(z):=
\ds\frac{n-1}{2(n+1)}z_1^2-z_2$ one gets the equality
$\ds\max_{\partial \Bbb G_n}|P_n|=M_n:=\frac{(n-1)(n+2)}{2(n+1)}.$
This means that if $$\ds g_n(t):=\frac{1}{2}\sum_{j=1}^n
t_j^2-\frac{1}{n+1}(\sum_{j=1}^n t_j)^2,\ t\in\Bbb C^n,$$ then
$\ds\max_{\Bbb T^n}|g_n|=M_n.$ Indeed, let $\ds
M_n^\ast=\max_{\Bbb T^n}|g_n|.$ Since
$g_n(e^{i\theta}t)=e^{2i\theta}g_n(t)$ for any $\theta\in\Bbb R,$
$t\in\Bbb C^2,$ there exists a point $u\in\Bbb T^n$ such that
$g_n(u)=M_n^\ast.$ Setting $u_j=x_j+iy_j,$ $x_j,y_j\in\Bbb R,$
$1\le j\le n,$ it follows that $$
M_n^\ast=\mbox{Re}(g_n(u))=\frac{1}{2}\sum_{j=1}^n
(x_j^2-y_j^2)+\frac{1}{n+1}((\sum_{j=1}^ny_j)^2-(\sum_{j=1}^nx_j)^2)$$
$$\le\frac{1}{2}\sum_{j=1}^n(x_j^2-y_j^2)+\frac{1}{n+1}(n\sum_{j=1}^ny_j^2-
(\sum_{j=1}^nx_j)^2)$$
$$=\frac{(n-1)n}{2(n+1)}+\frac{1}{n+1}(\sum_{j=1}^n
x_j^2-(\sum_{j=1}^nx_j)^2)$$ by the Cauchy-Schwarz inequality and
the equalities $y_1^2=1-x_1^2,\dots,y_n^2=1-x_n^2.$ The last term
is a linear function in any $x_j.$ Hence it attains maximum at
$\pm1.$ Since $n$ is odd, then $$
M_n^\ast=\frac{(n-1)n}{2(n+1)}+\frac{n-1}{n+1}=M_n$$ and the
maximum is attained at $t\in\Bbb T^n$ if and only $[n/2]$ or
$[n/2]+1$ of the $t_j$'s are equal to some $t_0\in\Bbb T$ and the
other ones to $-t_0.$

Using this last fact, it is not difficult to see that if
$\varepsilon>0$ is small and
$$g_{n,\varepsilon}(t)=g_n(t)+\varepsilon\sum_{j=1}^n
t_j^2-\varepsilon(n+1)(\sum_{j=1}^n t_j)^2,\ t\in\Bbb C^n,$$ then
$\ds\max_{\Bbb T^n}|g_{n,\varepsilon}|<M_n.$ Therefore, for
$$P_{n,\varepsilon}=\frac{n-1-2n(n+1)\varepsilon}{2(n+1)}z_1-(1+2\varepsilon)
z_2$$ one has the inequality $\ds\max_{\partial \Bbb
G_n}|P_{n,\varepsilon}|<M_n$ which implies that $$\ds\gamma_{\Bbb
G_n}(0;e_2)>\frac{1}{M_n}=\frac{2}{n}\left(1+\frac{2}{(n-1)(n+2)}\right).$$

{\it The upper bound:} In virtue of (1), we have to show that if
$c\in\Bbb C,$ then $$m_{n,c}:=\max_{z\in\partial
G_n}|z_2+cz_1^2|>\frac{n(n^2-1)}{2(n^2+1)}.$$ The coefficients of
the polynomials $(t-1)^n$ and $(t-1)(t^2-1)^{\frac{n-1}{2}}$ give
points $z\in\partial\Bbb G_n$ with $\ds z_1=n,\ds
z_2=\frac{n(n-1)}{2}$ and $\ds z_1=1,z_2=\frac{1-n}{2},$
respectively. Then $$2m_{n,c}\ge\max\{|n-1-2c|,|n(n-1)+2cn^2|\}$$
and hence $$2(n^2+1)m_{n,c}\ge|n^2(n-1)-2cn^2|+|n(n-1)+2cn^2|$$
$$\ge n^2(n-1)+n(n-1)=n(n^2-1).$$ This implies that $\ds
m_{n,c}\ge\frac{n(n^2-1)}{2(n^2+1)}.$ Assume that the equality
holds. Then $\ds c=-\frac{(n-1)^2}{2(n^2+1)}.$ On the other hand,
the coefficients of the polynomial $(t-i)(t-1)^{n-1}$ give a point
$z\in\partial\Bbb G_n$ with $\ds
z_1=n-1+i,z_2=\frac{(n-1)(n-2)}{2}+(n-1)i,$ for which
$\ds\left|z_2-\frac{(n-1)^2}{2(n^2+1)}z_1^2\right|>\frac{n(n^2-1)}{2(n^2+1)},$
a contradiction.
\end{proof}

\section{The proof of the inequality $\hat\gamma^{(2)}_{\Bbb G_3}(0;e_2)>\gamma_{\Bbb
G_3}(0;e_2)$}

Let $D$ be a domain in $\Bbb C^n$ and $k\in\Bbb N.$ Recall that
the $k$-th Reiffen pseudometric is defined as (see \cite{Jar-Pfl})
$$\gamma^{(k)}_D(z;X):=
\sup\{\left|\frac{f^{(k)}(z)X}{k!}\right|^{\frac{1}{k}}:f\in\mathcal
O(D,\Bbb D),\mbox{ord}f_z\ge k\}.$$ Note that
$\gamma_D\le\gamma_D^{(k)}\le\kappa_D.$ Denote by $\hat
\gamma^{(k)}_D$ the largest pseudonorm which does not exceed
$\gamma^{(k)}_D.$ Since $\gamma_D(z;\cdot)$ is a pseudonorm, it
follows that $\gamma_D\le\hat\gamma_D^{(k)}\le\hat\kappa_D.$ We
also point out that the family $O(\Bbb G_3,\Bbb D)$ is normal and
then the argument as in \cite{NP} shows that there are $m$ ($m\le
2n-1$) $\Bbb R$-linearly independent vectors $X_1,\dots,X_m\in\Bbb
C^n$ with the sum $X$ such that
$$\ds\hat\gamma^{(k)}_D(z;X)=\sum_{j=1}^m \gamma^{(k)}_D(z;X).$$

The purpose of this section is to show the following

\begin{proposition}
$\hat\gamma^{(2)}_{\Bbb G_3}(0;e_2)>\gamma_{\Bbb G_3}(0;e_2).$ In
particular, $\hat\kappa_{\Bbb G_3}(0;e_2)>\gamma_{\Bbb
G_3}(0;e_2)$ and hence $k_{\Bbb G_3}(0,\cdot)\neq c_{\Bbb
G_3}(0,\cdot).$
\end{proposition}

\noindent{\bf Remark.} We believe that the idea of the proof below
works for $\Bbb G_n$ for any $n\ge 3.$
\medskip

Proposition 5 is a consequence of the next two lemmas.

\begin{lemma} $\ds\gamma_{\Bbb G_3}(0;e_2)\le C_0:=\sqrt\frac{8}{13\sqrt
{13}-35}=0,8208\dots.$
\end{lemma}

\begin{lemma} $\ds\hat\gamma_{\Bbb G_3}^{(2)}(0;e_2)\ge
C_1=\sqrt{0,675}=0,8215\dots.$
\end{lemma}

\noindent{\it Proof of Lemma 6.} By (1), we have to show that for
any $c\in\Bbb C$ one has $$\max_{z\in\partial\Bbb
G_3}|z_2-cz_1^2|^2\ge\frac{1}{C_0^2}.$$ First, observe that it is
enough to prove this inequality in the case, when $c\in\Bbb R.$
Indeed, for any $z\in\partial\Bbb G_3$ one has that $\overline
z\in\partial\Bbb G_3$ and therefore $$2\max_{z\in\partial\Bbb
G_3}|z_2-cz_1^2|\ge\max_{z\in\partial\Bbb G_3}(|z_2-cz_1^2|+
|\overline z_2-c\overline z_1^2|)$$ $$\ge\max_{z\in\partial\Bbb
G_3}|2z_2-(c+\overline c)z_1^2|=2\max_{z\in\partial\Bbb
G_3}|z_2-\mbox{Re}(c)z_1^2|.$$

Let now $c\in\Bbb R.$ Then $$\max_{z\in\partial\Bbb
G_3}|z_2-cz_1^2|^2\ge\max_{\varphi\in[0,2\pi)}|1+2e^{i\varphi}-c(2+e^{i\varphi})^2|^2$$
$$=\max_{\varphi\in[0,2\pi)}(4c(4c-1)\cos^2\varphi+4(10c^2-7c+1)\cos\varphi+25c^2-22c+5).$$
Set $$f_c(x):=4c(4c-1)x^2+4(2c-1)(5c-1)x+25c^2-22c+5,\
x\in[-1,1].$$ If $\ds
c\not\in\Delta:=\left(\frac{1}{6},\frac{5-\sqrt{17}}{4}\right),$
then $$\max_{x\in[-1,1]}f_c(x)=\max\{f_c(-1),f_c(1)\}\ge
\left(\frac{9-\sqrt{17}}{4}\right)^2>\frac{1}{C_0^2}.$$ Otherwise,
$$\max_{x\in[-1,1]}f_c(x)=f_c\left(\frac{10c^2-7c+1}{2c(1-4c)}\right)
=\frac{(3c-1)^3}{c(4c-1)}=:g(c)$$ and it remains to check that
$\ds\min_{c\in\Delta}g(c)=g\left(\frac{\sqrt{13}-1}{12}\right)=\frac{1}{C_0^2}.$\qed
\smallskip

\noindent{\bf Remark.} Set $\ds c_0=\frac{\sqrt{13}-1}{12}$ and
$\ds M:=\max_{z\in\partial\Bbb G_3}|z_2-c_0z_1^2|.$ As in the
proof of Proposition 4 we have that $$M=\max_{z\in\partial\Bbb
G_3}\hbox{Re}(z_2-c_0z_1^2)=\max_{\alpha,\beta,\gamma\in\Bbb R}
h(\alpha,\beta,\gamma),$$ where
$$h(\alpha,\beta,\gamma)=(1-2c_0)(\cos(\alpha+\beta)+\cos(\beta+\gamma)
+\cos(\gamma+\alpha))$$
$$-c_0(\cos2\alpha+\cos2\beta+cos2\gamma).$$ Computer calculations
show that the critical points of $h$ (up to permutations of the
variables) are of the form ($k\pi,l\pi,m\pi)$ or
$(\pm\alpha_0+j\pi/2+2k\pi,\pm\alpha_0+j\pi/2+2l\pi,\pm\gamma_0+j\pi/2+2m\pi),$
$k,l,m\in\Bbb Z,$ $j=0,1,2,3.$ Then it follows by the proof of
Lemma 6 that $\ds M=C_0^{-1}$ which implies that in fact
$\ds\gamma_{\Bbb G_3}(0;e_2)=C_0.$
\medskip

\noindent{\it Proof of Lemma 7.} Let
$$f(z)=0,675z_2^2-0,291z_2z_1^2+0,033z_1^4.$$  We claim that
$\ds\max_{z\in\partial\Bbb G_3}|f(z)|<1.$ Set
$\theta=(\theta_1,\theta_2),\ \theta_1,\theta_2\in[0,2\pi),$
$$g_1(\theta)=1+e^{i\theta_1}+e^{i\theta_2},\
g_2(\theta)=e^{i(\theta_1+\theta_2)}+e^{i\theta_1}+e^{i\theta_2},$$
$$g(\theta)=0,675g_2^2(\theta)-0,291
g_2(\theta)g_1^2(\theta)+0,033g_1^4(\theta).$$ We have to show
that $\max|g(\theta)|<1.$ Set
$$d(\theta,\tilde\theta)=\max\{|\theta_1-\tilde\theta_1|,
|\theta_2-\tilde\theta_2|\}.$$ Since
$|e^{i\theta_j}-e^{i\tilde\theta_j}|\le|\theta_j-\tilde\theta_j|,$
$j=1,2,$ then
$$|g_1(\theta)-g_1(\tilde\theta)|\le2d(\theta,\tilde\theta),\
|g_2(\theta)-g_2(\tilde\theta)|\le4d(\theta,\tilde\theta).$$ Now
the inequalities $|g_1|\le 3,|g_2|\le 3$ imply that
$$|g(\theta)-g(\tilde\theta)|\le(0,675\cdot 24+0,291\cdot 72+0,033\cdot 216)
d(\theta,\tilde\theta)=44,28d(\theta,\tilde\theta).$$

Let now $\theta_1,\theta_2$ vary on the interval
$[0;6,2832]\supset[0,2\pi]$ with step $4\cdot 10^{-5}.$ A simple
computer programme shows that $|g(\theta)|\le 0,999$ for all
running $\theta=(\theta_1,\theta_2).$ (In fact, one may conjecture
that $\max |g(\theta)|=0,999$ and this maximum is attained at the
points $(0,\pi),$ $(\pi,0)$ and $(\pi,\pi).$) Then the
inequalities $|g(\theta)-g(\tilde\theta)|\le
44,28d(\theta,\tilde\theta)$ and $\ds\frac{2}{44,28}\cdot
10^{-3}>4\cdot 10^{-5}$ easily prove that $\max |g(\theta)|<1.$

It follows that if $X\in\Bbb C^3$ is in the span of $e_1$ and
$e_3,$ then $f$ is a competitor for $\gamma^{(2)}_{\Bbb
G_3}(0;e_2+X)$ and hence $\gamma^{(2)}_{\Bbb G_3}(0;e_2+X)\ge
C_1.$

On the other hand, recall that we may find five vectors
$X_1,\dots,X_5\in\Bbb C^3$ (possible some of them 0) with sum
$e_2$ such that $\ds\hat\gamma^{(2)}_{\Bbb
G_3}(0;e_2)=\sum_{j=1}^5\gamma^{(2)}_{\Bbb G_3}(0;X_j).$ Since
$\gamma^{(2)}_{\Bbb G_3}(0;X_j)\ge|(e_2,\overline{X_j})|C_1,$ then
$\ds\hat\gamma^{(2)}_{\Bbb G_3}(0;e_2)\ge C_1.$\qed

Finally, we point out that $\gamma^{(2)}_{\Bbb G_n}(0;\cdot)$ is
not a norm.

\begin{proposition} If $X_1,X_n\in\Bbb C,$ then $$\gamma^{(2)}_{\Bbb
G_n}(0;X_1e_1+X_ne_n)\ge\sqrt{\frac{n+1}{2}\gamma_{\Bbb
G_n}(0;e_2)|X_1X_n|}.$$

In particular, since $\ds\gamma_{\Bbb G_3}(0;e_2)>\frac{2}{3}$ and
$\ds\gamma_{\Bbb G_n}(0;e_n)\ge\frac{2}{n},$ then
$$\gamma^{(2)}_{\Bbb G_n}(0;ne_1+e_n)>2=\hat\kappa_{\Bbb
G_n}(0;ne_1+e_n)=\gamma^{(2)}_{\Bbb
G_n}(0;ne_1)+\gamma^{(2)}_{\Bbb G_n}(0;e_n),\ n\ge 3.$$
\end{proposition}

\begin{proof} Let $t_1,\dots,t_n\in\Bbb D.$ Consider
$\ds\sum_{k=1}^n\frac{t_k^{n+1}}{n}$ as a function $f_n$ of
$z_1,\dots,z_n.$ Then $f_n\in\mathcal O(\Bbb G_n,\Bbb D),$
$\mbox{ord}_0 f_n=2$ and by the Waring formula (cf. \cite{Waer}))
the coefficient at $z_1z_n$ equals $\ds(-1)^{n-1}\frac{n+1}{n}.$
Hence $$\gamma^{(2)}_{\Bbb
G_n}(0;X_1e_1+X_ne_n)\ge\sqrt{\frac{n+1}{n}\gamma_{\Bbb
G_n}(0;e_2)|X_1X_n|}.$$ Since $\ds\gamma_{\Bbb
G_n}(0;e_2)=\frac{2}{n}$ for even $n,$ we are done for such $n.$

On the other hand, we know by Proposition 3 that there is $c_n$
such that $P$ with $P(z):=2C_nz_2-c_nz_1^2$ is an extremal
function for $\gamma_{\Bbb G_n}(0;e_2)=:2C_n.$ For $n=2k-1$
replace $t_1,\dots,t_n$ by $t_1^k,\dots,t_n^k.$ Then we obtain the
function $$g_n(\sigma_n(t)):=\tilde
g_n(t)=(C_n-c_n)\left(\sum_{j=1}^nt_j^k\right)^2
-C_n\sum_{j=1}^nt_j^{2k}.$$ Then $g_n\in\mathcal O(\Bbb G_n,\Bbb
D),$ $\mbox{ord}_0 g_n=2,$ and the coefficient at $z_1z_n$ equals
$-(n+1)C_n.$ Now, it is enough to take $g_n$ as a competitor for
$\gamma^{(2)}_{\Bbb G_{n}}(0;X_1e_1+X_ne_n).$
\end{proof}

\end{document}